\documentclass{amsart}

\newtheorem{theorem}{Theorem}[section]
\newtheorem{lemma}[theorem]{Lemma}

\newtheorem{proposition}[theorem]{Proposition}
\newtheorem{remark}[theorem]{Remark}

\theoremstyle{definition}

\theoremstyle{remark}

\numberwithin{equation}{section}

\begin{document}

\title[Power Moments of Kloosterman Sums  with Square Arguments]{    $\begin{array}{c}
  \text{Ternary Codes Associated with $O^{-}(2n,q)$ and }\\
 \text{ Power Moments of Kloosterman Sums  with Square Arguments} \\
\end{array} $            }

\author{dae san kim}
\address{Department of Mathematics, Sogang University, Seoul 121-742, Korea}
\curraddr{Department of Mathematics, Sogang University, Seoul
121-742, Korea}
\email{dskim@sogong.ac.kr}
\thanks{}

%    General info
\subjclass[]{}

\date{}

\dedicatory{ }


\begin{abstract}
In this paper, we construct three ternary linear codes associated
with the orthogonal group $O^{-}(2,q)$ and the special orthogonal
groups $SO^- (2,q)$ and $SO^{-}(4,q)$. Here $q$ is a power of
three. Then we obtain recursive formulas for the power moments of
Kloosterman sums with square arguments and for the even power
moments of those in terms of the frequencies of weights in the
codes. This is done via Pless power moment identity and by
utilizing the explicit expressions of ``Gauss sums" for the
orthogonal and special orthogonal groups
$O^{-}(2n,q)$ and $SO^{-}(2n,q)$.\\

Key words-ternary linear code, power moment, Kloosterman sum,
square argument, Pless power moment identity, Gauss sum,
orthogonal group, weight distribution.\\

\end{abstract}

\maketitle

%%%%%%%%%%%%%%%%%%%%%%%%%%%%%%%%%%%%%%%%%%%%%%%%%%%%%%%%%%%%%%%%%%%%%%%%
%%%%%%%%%%%%%%%%%%%%%%%%%%%%%%%%%%%%%%%%%%%%%%%%%%%%%%%%%%%%%%%%%%%%%%%%
\section{Introduction}
%%%%%%%%%%%%%%%%%%%%%%%%%%%%%%%%%%%%%%%%%%%%%%%%%%%%%%%%%%%%%%%%%%%%%%%%
%%%%%%%%%%%%%%%%%%%%%%%%%%%%%%%%%%%%%%%%%%%%%%%%%%%%%%%%%%%%%%%%%%%%%%%%

Let $\psi$ be a nontrivial additive character of the finite field
$\mathbb{F}_{q}$ with $q=p^r$ elements ($p$ a prime). Then the
Kloosterman sum $K(\psi;a)$ (\cite{LN1}) is defined by

\begin{align*}
K(\psi;a)=\sum_{\alpha\in\mathbb{F}_{q}^{*}}\psi(\alpha_{}^{}+a\alpha_{}^{-1})
(a\in\mathbb{F}_{q}^{*}).
\end{align*}
For this, we have the Weil bound
%(1)%%%%%%%%%%%%%%%%%%%%%%%%%%%%%%%%%%%%%%%%%%%%%%%%%%%%%%%%%%%%%%
\begin{equation}\label{a}
|K(\psi;a)|\leq 2\sqrt{q}.
\end{equation}
%%%%%%%%%%%%%%%%%%%%%%%%%%%%%%%%%%%%%%%%%%%%%%%%%%%%%%%%%%%%%%%%%%
The Kloosterman sum was introduced in 1926(\cite{K1}) to give an
estimate for the Fourier coefficients of modular forms.

For each nonnegative integer $h$, by $MK(\psi)^{h}$ we will denote
the $h$-th moment of the Kloosterman sum $K(\psi;a)$. Namely, it is
given by
\begin{align*}
MK(\psi)^{h}=\sum_{a\in\mathbb{F}_{q}^{*}}K(\psi;a)^{h}.
\end{align*}
If $\psi=\lambda$ is the canonical additive character of
$\mathbb{F}_{q}$, then $MK(\lambda)^{h}$ will be simply denoted by
$MK^{h}$.

Also, we introduce an incomplete power moments of Kloosterman sums.
Namely, for every nonnegative integer $h$, and $\psi$ as before, we
define
%2%%%%%%%%%%%%%%%%%%%%%%%%%%%%%%%%%%%%%%%%%%%%%%%%%%%%%%%%%%%%%%%%
\begin{equation}\label{b}
SK(\psi)^{h}=\sum_{a\in\mathbb{F}_{q}^{*},~a~square}K(\psi;a)^{h},
\end{equation}
%%%%%%%%%%%%%%%%%%%%%%%%%%%%%%%%%%%%%%%%%%%%%%%%%%%%%%%%%%%%%%%%%%
which is called the $h$-th moment of Kloosterman sums with
``square arguments". If $\psi=\lambda$ is the canonical additive
character of $\mathbb{F}_q$, then $SK(\lambda)^{h}$ will be
denoted by $SK^h$, for brevity.

Explicit computations on power moments of Kloosterman sums were
begun with the paper \cite{S1} of Sali\'{e} in 1931, where he
showed, for any odd prime $q$,
\begin{align*}
MK^h=q^2 M_{h-1} -(q-1)^{h-1}+2(-1)^{h-1}~(h\geq1).
\end{align*}
Here $M_0=0$, and for $h\in\mathbb{Z}_{>0}$,
\begin{align*}
M_h=|\{(\alpha_1,\cdots,\alpha_h)\in(\mathbb{F}_q^*)^h|\sum_{j=1}^h
\alpha_j=1=\sum_{j=1}^h \alpha_j^{-1}\}|.
\end{align*}
For $q=p$ odd prime, Sali\'{e} obtained $MK^1$, $MK^2$, $MK^3$,
$MK^4$ in \cite{S1} by determining $M_1$, $M_2$, $M_3$. On the
other hand, $MK^5$ can be expressed in terms of the $p$-th
eigenvalue for a weight 3 newform on $\Gamma_0(15)$(cf. \cite{L1},
\cite{PV1}). $MK^6$ can be expressed in terms of the $p$-th
eigenvalue for a weight 4 newform on $\Gamma_0(6)$(cf.
\cite{HS1}). Also, based on numerical evidence, in \cite{E1} Evans
was led to propose a conjecture which expresses $MK^7$ in terms of
Hecke eigenvalues for a weight 3 newform on $\Gamma_0(525)$ with
quartic nebentypus of conductor 105.

Assume now that $q=3^r$. Recently, Moisio was able to find
explicit expressions of $MK^h$, for $h\leq10$(cf.\cite{M1}). This
was done, via Pless power moment identity, by connecting moments
of Kloosterman sums and the frequencies of weights in the ternary
Melas code of length $q-1$, which were known by the work of Geer,
Schoof and Vlugt in \cite{GS1}. In this paper, we will be able to
produce three recursive formulas generating power moments of
Kloosterman sums with square arguments over finite fields of
characteristic three. To do that, we will construct three ternary
linear codes $C(SO^- (2,q))$, $C(O^- (2,q))$, and $C(SO^- (4,q))$,
respectively associated with $SO^- (2,q)$, $O^- (2,q)$, and $SO^-
(4,q)$, and express those power moments in terms of the
frequencies of weights in each code. Then, thanks to our previous
results on the explicit expressions of ``Gauss sums" for the
orthogonal and special orthogonal groups $O^- (2n,q)$ and $SO^-
(2n,q)$, we can express the weight of each codeword in the duals
of the codes in terms of Kloosterman sums with square arguments.
Then our formulas will follow immediately from the Pless power
moment identity (cf. (\ref{q1})). Similar results of this paper
were obtained in \cite{DJ1} for the case of finite symplectic
groups over finite fields of characteristic three. Also, in the
same case infinite families of recursive formulas were derived in
\cite{D5} by using explicit expressions of exponential sums
associated with certain double cosets.

Theorem \ref{Aa} in the following (cf. (\ref{e}), (\ref{f}),
(\ref{h})-(\ref{l})) is the main result of this paper. Henceforth,
we agree that, for nonnegative integers $a$, $b$, $c$,
%(3)%%%%%%%%%%%%%%%%%%%%%%%%%%%%%%%%%%%%%%%%%%%%%%%%%%%%%%%%%%%%%%%%%%%
\begin{equation}\label{c}
{\binom{c}{a,b}}={\frac{c!}{a!~b!~(c-a-b)!}},~if~a+b\leq c,
\end{equation}
%%%%%%%%%%%%%%%%%%%%%%%%%%%%%%%%%%%%%%%%%%%%%%%%%%%%%%%%%%%%%%%%%%%%%%%
and
%(4)%%%%%%%%%%%%%%%%%%%%%%%%%%%%%%%%%%%%%%%%%%%%%%%%%%%%%%%%%%%%%%%%%%%
\begin{equation}\label{d}
{\binom{c}{a,b}}=0,~if~a+b>c.
\end{equation}
%%%%%%%%%%%%%%%%%%%%%%%%%%%%%%%%%%%%%%%%%%%%%%%%%%%%%%%%%%%%%%%%%%%%%%%
\begin{theorem}\label{Aa}
Let $q=3^r$. Then we have the following. \\$(1)$ For
$h$=1,2,$\cdots$,
\begin{equation}\label{e}
\begin{split}&SK_{}^{h}=-\sum_{j=0}^{h-1}{\binom{h}{j}}(q+1)^{h-j}SK_{}^{j}\\
            &~\quad\quad +q\sum_{j=0}^{min\{N_1,h\}}(-1)^{j}C_{1,j}
            \sum_{t=j}^{h}t!S(h,t)3^{h-t}2^{t-h-j-1}{\binom{N_1 -j}{N_1
            -t}},
\end{split}
\end{equation}
where $N_1=|SO^- (2,q)|=q+1$, and $\{C_{1,j}\}_{j=0}^{N_1}$ is the
weight distribution of the ternary linear code $C(SO^- (2,q))$
given by
\begin{equation}\label{f}
C_{1,j}=\sum{\binom{1}{\nu_1,\mu_1}}{\binom{1}{\nu_{-1},\mu_{-1}}}
\prod_{\beta^2-1~nonsquare}{\binom{2}{\nu_{\beta},\mu_{\beta}}}
~(j=0,\cdots,N_1).
\end{equation}
Here the sum is over all the sets of nonnegative integers
$\{\nu_{\beta}\}_{\beta^2-1~nonsquare}\cup\{\nu_{\pm1}\}$ and
$\{\mu_{\beta}\}_{\beta^2-1~nonsquare}\cup\{\mu_{\pm1}\}$
satisfying
\begin{align*}
\nu_1+\nu_{-1}+\sum_{\beta^2-1~nonsquare}\nu_{\beta}+\mu_1+\mu_{-1}
+\sum_{\beta^2-1~nonsquare}\mu_{\beta}=j
\end{align*}
and
\begin{align*}
\nu_1
1+\nu_{-1}(-1)+\sum_{\beta^2-1~nonsquare}\nu_{\beta}\beta=\mu_1
1+\mu_{-1}(-1)+\sum_{\beta^2-1~nonsquare}\mu_{\beta}\beta.
\end{align*}
In addition, $S(h,t)$ is the Stirling number of the second kind
defined by
%(7)%%%%%%%%%%%%%%%%%%%%%%%%%%%%%%%%%%%%%%%%%%%%%%%%%%%%%%%%%%%%%%%%%%%%%
\begin{equation}\label{g}
S(h,t)=\frac{1}{t!}\sum_{j=0}^{t}(-1)^{t-j}{\binom{t}{j}}j^h.
\end{equation}\\$(2)$ For $h=1,2,\cdots,$
%(8)%%%%%%%%%%%%%%%%%%%%%%%%%%%%%%%%%%%%%%%%%%%%%%%%%%%%%%%%%%%%%%%%%%%
\begin{equation}\label{h}
\begin{split}&SK_{}^{h}=-\sum_{j=0}^{h-1}{\binom{h}{j}}(q+1)^{h-j}SK_{}^{j}\\
            &~\quad\quad +q\sum_{j=0}^{min\{N_2,h\}}(-1)^{j}C_{2,j}
            \sum_{t=j}^{h}t!S(h,t)3^{h-t}2^{t-h-j-1}{\binom{N_2
            -j}{N_2-t}},
\end{split}
\end{equation}
%%%%%%%%%%%%%%%%%%%%%%%%%%%%%%%%%%%%%%%%%%%%%%%%%%%%%%%%%%%%%%%%%%%%%%%
where $N_2=|O^- (2,q)|=2(q+1)$, and $\{C_{2,j}\}_{j=0}^{N_2}$ is
the weight distribution of the ternary linear code $C(O^- (2,q))$
given by:\\$(a)$ For $r$ even,
\begin{equation}\label{i}
\begin{split}
&C_{2,j}=\sum{\binom{q+1}{\nu_0,\mu_0}}{\binom{1}{\nu_1,\mu_1}}{\binom{1}{\nu_{-1},\mu_{-1}}}
\prod_{\beta^2-1~nonsquare}{\binom{2}{\nu_{\beta},\mu_{\beta}}}~(j=0,\cdots,N_2).
\end{split}
\end{equation}
 Here the sum is over all the sets of nonnegative integers
$\{\nu_{\beta}\}_{\beta^2-1~nonsquare}$\\$\cup\{\nu_{\pm1}\}\cup\{\nu_0\}$
and
$\{\mu_{\beta}\}_{\beta^2-1~nonsquare}\cup\{\mu_{\pm1}\}\cup\{\mu_0\}$
satisfying
\begin{align*}
\begin{split}
&\nu_0+\nu_1+\nu_{-1}+\sum_{\beta^2-1~nonsquare}\nu_{\beta}\\
&\qquad\qquad\qquad+\mu_0+\mu_1+\mu_{-1}
+\sum_{\beta^2-1~nonsquare}\mu_{\beta}=j
\end{split}
\end{align*}
and
\begin{align*}
\begin{split}
&\nu_1 1+\nu_{-1}(-1)+\sum_{\beta^2-1~nonsquare}\nu_{\beta}\beta\\
&\qquad\qquad\qquad=\mu_1
1+\mu_{-1}(-1)+\sum_{\beta^2-1~nonsquare}\mu_{\beta}\beta.
\end{split}
\end{align*}
\\$(b)$   For $r$ odd,
%(10)%%%%%%%%%%%%%%%%%%%%%%%%%%%%%%%%%%%%%%%%%%%%%%%%%%%%%%%%%%%%%%%%%%%%
\begin{equation}\label{j}
\begin{split}
C_{2,j}=\sum{\binom{1}{\nu_1,\mu_1}}{\binom{1}{\nu_{-1},\mu_{-1}}}
\prod_{\beta^2-1\neq-1~nonsquare}{\binom{2}{\nu_{\beta},\mu_{\beta}}}{\binom{q+3}{\nu_0,\mu_0}}
\end{split}
\end{equation}
~~~~~~$ (j=0,\cdots,N_2).$\\\\
%%%%%%%%%%%%%%%%%%%%%%%%%%%%%%%%%%%%%%%%%%%%%%%%%%%%%%%%%%%%%%%%%%%%%%%%
Here the sum is over all the sets of nonnegative integers
$\{\nu_{\beta}\}_{\beta^2-1\neq-1~nonsquare}$\\$\cup\{\nu_{\pm1}\}\cup\{\nu_0\}$
and
$\{\mu_{\beta}\}_{\beta^2-1\neq-1~nonsquare}\cup\{\mu_{\pm1}\}\cup\{\mu_0\}$
satisfying
\begin{align*}
\begin{split}
&\nu_0+\nu_1+\nu_{-1}+\sum_{\beta^2-1\neq-1~nonsquare}\nu_{\beta}\\
&\qquad\qquad\qquad+\mu_0+\mu_1+\mu_{-1}
+\sum_{\beta^2-1\neq-1~nonsquare}\mu_{\beta}=j
\end{split}
\end{align*}
and
\begin{align*}
\begin{split}
&\nu_1 1+\nu_{-1}(-1)+\sum_{\beta^2-1\neq-1~nonsquare}\nu_{\beta}\beta\\
&\qquad\qquad\qquad=\mu_1
1+\mu_{-1}(-1)+\sum_{\beta^2-1\neq-1~nonsquare}\mu_{\beta}\beta.
\end{split}
\end{align*}
\\$(3)$
For $h=1,2,\cdots,$
%(11)%%%%%%%%%%%%%%%%%%%%%%%%%%%%%%%%%%%%%%%%%%%%%%%%%%%%%%%%%%%%%%%%%%%%
\begin{equation}\label{k}
\begin{split}&SK_{}^{2h}=-\sum_{j=0}^{h-1}{\binom{h}{j}}(q^4+q^3-q-1)^{h-j}SK_{}^{2j}\\
            &~\quad\quad +q^{1-2h}\sum_{j=0}^{min\{N_3,h\}}(-1)^{j}C_{3,j}
            \sum_{t=j}^{h}t!S(h,t)3^{h-t}2^{t-h-j-1}{\binom{N_3
            -j}{N_3 -t}},
\end{split}
\end{equation}
%%%%%%%%%%%%%%%%%%%%%%%%%%%%%%%%%%%%%%%%%%%%%%%%%%%%%%%%%%%%%%%%%%%%%%%%
where $N_3=|SO^- (4,q)|=q^2 (q^4 -1)$, and
$\{C_{3,j}\}_{j=0}^{N_3}$ is the weight distribution \\of the
ternary linear code $C(SO^- (4,q))$ given by

%(12)%%%%%%%%%%%%%%%%%%%%%%%%%%%%%%%%%%%%%%%%%%%%%%%%%%%%%%%%%%%%%%%%%%%%
\begin{equation}\label{l}
\begin{split}
&\qquad C_{3,j}=\sum{\binom{-q^2 \delta(2,q;0)+q^4 +2q^3
-3q^2}{\nu_0,\mu_0}}\\
&\qquad\qquad\qquad\qquad\qquad\times\prod_{\beta\in\mathbb{F}_{q}^{*}}{\binom{-q^2
\delta(2,q;\beta)+q^5 +q^4 +q^3 -3q^2}
{\nu_{\beta},\mu_{\beta}}}\\
&\qquad\qquad\qquad\qquad\qquad\qquad\qquad\qquad\qquad\qquad\qquad\qquad(j=0,\cdots,N_3).
\end{split}
\end{equation}
%%%%%%%%%%%%%%%%%%%%%%%%%%%%%%%%%%%%%%%%%%%%%%%%%%%%%%%%%%%%%%%%%%%%%%%%
Here the sum runs over all the sets of nonnegative integers
$\{\nu_\beta\}_{\beta\in\mathbb{F}_q}$ and
$\{\mu_\beta\}_{\beta\in\mathbb{F}_q}$ satisfying
$\sum_{\beta\in\mathbb{F}_q}\nu_\beta+\sum_{\beta\in\mathbb{F}_q}\mu_\beta=j$
and $\sum_{\beta\in\mathbb{F}_q}\nu_\beta
\beta=\sum_{\beta\in\mathbb{F}_q}\mu_\beta \beta$, and, for every
$\beta\in\mathbb{F}_{q}$,
$\delta(2,q;\beta)=|\{(\alpha_1,\alpha_2)\in(\mathbb{F}_{q}^{*})^2|\alpha_1+\alpha_{1}^{-1}+
\alpha_{2}^{}+\alpha_{2}^{-1}=\beta\}|$.
\end{theorem}

%Theorem1%%%%%%%%%%%%%%%%%%%%%%%%%%%%%%%%%%%%%%%%%%%%%%%%%%%%%%%%%%%%%%

%%%%%%%%%%%%%%%%%%%%%%%%%%%%%%%%%%%%%%%%%%%%%%%%%%%%%%%%%%%%%%%%%%%%%%%%
%%%%%%%%%%%%%%%%%%%%%%%%%%%%%%%%%%%%%%%%%%%%%%%%%%%%%%%%%%%%%%%%%%%%%%%%
\section{$O^{-}(2n,q)$}
%%%%%%%%%%%%%%%%%%%%%%%%%%%%%%%%%%%%%%%%%%%%%%%%%%%%%%%%%%%%%%%%%%%%%%%%
%%%%%%%%%%%%%%%%%%%%%%%%%%%%%%%%%%%%%%%%%%%%%%%%%%%%%%%%%%%%%%%%%%%%%%%%
For more details about this section, one is referred to the paper
\cite{D1} and \cite{D2}. Throughout this paper, the following
notations will be used:
\begin{align*}
\begin{split}
q&=3^r~(r\in\mathbb{Z}_{>0}),\\
\mathbb{F}_q&=~the~finite~field~ with~ q~ elements,\\
TrA&=~the~ trace~ of~ A~ for~ a~ square~ matrix~ A,\\
^t B&=~the~transpose~ of~ B~for~any~matrix~B.
\end{split}
\end{align*}

The orthogonal group $O^-(2n,q)$ over the field $\mathbb{F}_q$ is
defined as:
\begin{align*}
O^-(2n,q)=\{w\in GL(2n,q)|^twJw=J\},
\end{align*}
where
\begin{align*}
J=\left[
  \begin{array}{cccc}
    0 & 1_{n-1} & 0 & 0 \\
    1_{n-1} & 0 & 0 & 0 \\
    0 & 0 & 1 & 0 \\
    0 & 0 & 0 & -\epsilon \\
  \end{array}
\right],
\end{align*}
and $\epsilon$ is a fixed element in $\mathbb{F}_q^*\setminus
{\mathbb{F}_q^*}^2$, here and throughout this paper. \\
For convenience, we put
%(13)%%%%%%%%%%%%%%%%%%%%%%%%%%%%%%%%%%%%%%%%%%%%%%%%%%%%%%%%%%%%
\begin{equation}\label{m}
\delta_\epsilon=\left[
                  \begin{array}{cc}
                    1 & 0 \\
                    0 & -\epsilon \\
                  \end{array}
                \right].
 \end{equation}
%%%%%%%%%%%%%%%%%%%%%%%%%%%%%%%%%%%%%%%%%%%%%%%%%%%%%%%%%%%%%%%%%%
Then $O^-(2n,q)$ consists of all matrices

\begin{align*}
\left[
  \begin{array}{ccc}
    A & B & e \\
    C & D & f \\
    g & h & i \\
  \end{array}
\right](A,B,C,D~(n-1)\times(n-1),e,f~(n-1)\times2,g,h~2\times(n-1),i~2\times2)
\end{align*}
in $GL(2n,q)$ satisfying the relations:

\begin{align*}
\begin{split}
&{}^tAC+{}^tCA + {}^tg\delta_\epsilon g=0,\\
&{ }^tBD+{}^tDB +{ }^th\delta_\epsilon h=0,\\
&{ }^tef +{ }^tfe +{}^ti\delta_\epsilon i=\delta_\epsilon,\\
&{ }^tAD +{ }^tCB +{}^tg\delta_\epsilon h=1_{n-1},\\
&{ }^tAf +{ }^tCe +{}^tg\delta_\epsilon i=0,\\
&{ }^tBf +{ }^tDe +{}^th\delta_\epsilon i=0.\\
\end{split}
\end{align*}

 The special orthogonal group $SO^-(2n,q)$ over the field
$\mathbb{F}_q$ is defined as:
\begin{align*}
SO^-(2n,q)=\{w\in O^-(2n,q)|det~ w=1\},
\end{align*}
which is a subgroup of index 2 in $O^-(2n,q)$.

In particular, we have
\begin{align*}
O^-(2,q)&=\{i\in GL(2,q)|^ti\delta_\epsilon i=\delta_\epsilon\}
\end{align*}
%(14)%%%%%%%%%%%%%%%%%%%%%%%%%%%%%%%%%%%%%%%%%%%%%%%%%%%%%%%%%%%%%
\begin{equation}\label{n}
\qquad \qquad \qquad \qquad \qquad =SO^-(2,q) \amalg\left[
                   \begin{array}{cc}
                     1 & 0 \\
                     0 & -1 \\
                   \end{array}
                 \right]SO^-(2,q),~
\end{equation}
%%%%%%%%%%%%%%%%%%%%%%%%%%%%%%%%%%%%%%%%%%%%%%%%%%%%%%%%%%%%%%%%%%
with
\begin{align*}
\begin{split}
SO^-(2,q)&=\left\{ \left[
               \begin{array}{cc}
                 a & b\epsilon \\
                 b & a \\
               \end{array}
             \right]
\Big{|}~a,b\in\mathbb{F}_q,~a^2-b^2 \epsilon =1\right\}\\
&=\left\{\left[
      \begin{array}{cc}
        a & b\epsilon \\
        b & a \\
      \end{array}
    \right]
\Big{|}
~a+b\epsilon\in\mathbb{F}_q(\epsilon)~with~N_{\mathbb{F}_q(\epsilon)/\mathbb{F}_q}(a+b\epsilon)=1\right\}.
\end{split}
\end{align*}

Let $P(2n,q)$ be the maximal parabolic subgroup of $O^-(2n,q)$ given
by
\begin{align*}
\begin{split}
P&=P(2n,q)\\
&=\left\{\left[
      \begin{array}{ccc}
        A & 0 & 0 \\
        0 & ^tA^{-1} & 0 \\
        0 & 0 & i \\
      \end{array}
    \right]\left[
             \begin{array}{ccc}
               1_{n-1} & B & -^th\delta_\epsilon \\
               0 & 1_{n-1} & 0 \\
               0 & h & 1_2 \\
             \end{array}
           \right]
\Bigg{|}\begin{array}{c}
   A\in GL(n-1,q) \\
   i\in O^-(2n,q) \\
   ^tB+B+^th\delta_\epsilon h=0
 \end{array}
\right\},
\end{split}
\end{align*}
and let $Q=Q(2n,q)$ be the subgroup of $P(2n,q)$ of index 2
defined by
\begin{align*}
\begin{split}
Q&=Q(2n,q)\\
&=\left\{\left[
      \begin{array}{ccc}
        A & 0 & 0 \\
        0 & ^tA^{-1} & 0 \\
        0 & 0 & i \\
      \end{array}
    \right]\left[
             \begin{array}{ccc}
               1_{n-1} & B & -^th\delta_\epsilon \\
               0 & 1_{n-1} & 0 \\
               0 & h & 1_2 \\
             \end{array}
           \right]
\Bigg{|}\begin{array}{c}
   A\in GL(n-1,q) \\
   i\in SO^-(2n,q) \\
   ^tB+B+^th\delta_\epsilon h=0
 \end{array}
\right\}.
\end{split}
\end{align*}

From (\ref{n}), we see that
%(15)%%%%%%%%%%%%%%%%%%%%%%%%%%%%%%%%%%%%%%%%%%%%%%%%%%%%%%%%%%%%%%%%%%%
\begin{equation}\label{o}
P=Q\amalg\rho Q,
\end{equation}
%%%%%%%%%%%%%%%%%%%%%%%%%%%%%%%%%%%%%%%%%%%%%%%%%%%%%%%%%%%%%%%%%%%%%%%%
with
\begin{align*}
\rho=\left[
       \begin{array}{cccc}
         1_{n-1} & 0 & 0 & 0 \\
         0 & 1_{n-1} & 0 & 0 \\
         0 & 0 & 1 & 0 \\
         0 & 0 & 0 & -1 \\
       \end{array}
     \right]
.
\end{align*}

The Bruhat decomposition of $O^-(2n,q)$ with respect to $P=P(2n,q)$
is given by
\begin{align*}
O^-(2n,q)=\coprod_{r=0}^{n-1}P\sigma_r
P=\coprod_{r=0}^{n-1}P\sigma_r Q,
\end{align*}
which can be modified to give
%(16)%%%%%%%%%%%%%%%%%%%%%%%%%%%%%%%%%%%%%%%%%%%%%%%%%%%%%%%%%%%%%%%%%%%
\begin{equation}\label{p}
O^-(2n,q)=\coprod_{r=0}^{n-1}P\sigma_r (B_r\setminus Q),
\end{equation}
%%%%%%%%%%%%%%%%%%%%%%%%%%%%%%%%%%%%%%%%%%%%%%%%%%%%%%%%%%%%%%%%%%%%%%%%%
with
\begin{align*}
B_r=B_r (q)=\{w\in Q(2n,q)|\sigma_r w \sigma_r^{-1}\in P\}.
\end{align*}
Here $\sigma_r$ denotes the following matrix in $O^-(2n,q)$
\begin{align*}
\sigma_r=\left[
           \begin{array}{ccccc}
             0 & 0 & 1_r & 0 & 0 \\
             0 & 1_{n-1-r} & 0 & 0 & 0 \\
             1_r & 0 & 0 & 0 & 0 \\
             0 & 0 & 0 & 1_{n-1-r} & 0 \\
             0 & 0 & 0 & 0 & 1_2 \\
           \end{array}
         \right]
~(0\leq r\leq n-1).
\end{align*}
  One can also show that
\begin{align*}
|B_r\setminus Q|=   [{{}_{~~~r}^{n-1}}]_q q^{r(r+3)/2}(0\leq r\leq
n-1)~~\qquad~(cf. ~[4] ,~ (3.12), ~(3.21)),
\end{align*}
where, for integers $n$, $r$ with $0\leq r\leq n$, the
$q$-binomial coefficients are defined as:
\begin{align*}
 \Big[{{}_{~r}^{~n}
  \Big]}_q=\prod_{j=0}^{r-1}(q^{n-j}-1)/(q^{r-j}-1).
\end{align*}

Taking the decomposition in (\ref{o}) into consideration, we see
that (\ref{p}) can further be modified as
%(17)%%%%%%%%%%%%%%%%%%%%%%%%%%%%%%%%%%%%%%%%%%%%%%%%%%%%%%%%%%%%%%
\begin{equation}\label{q}
O^-(2n,q)=\coprod_{r=0}^{n-1}Q\sigma_r(B_r\setminus Q)
\amalg\coprod_{r=0}^{n-1}(\rho Q)\sigma_r(B_r\setminus Q),
\end{equation}
%%%%%%%%%%%%%%%%%%%%%%%%%%%%%%%%%%%%%%%%%%%%%%%%%%%%%%%%%%%%%%%%%
and
%(18)%%%%%%%%%%%%%%%%%%%%%%%%%%%%%%%%%%%%%%%%%%%%%%%%%%%%%%%%%%%%%%
\begin{equation}\label{r}
SO^-(2n,q)=\coprod_{\substack{0\leq r\leq
n-1\\r:even}}Q\sigma_r(B_r\setminus Q)
\amalg\coprod_{\substack{0\leq r\leq n-1\\r:odd}}(\rho
Q)\sigma_r(B_r\setminus Q).
\end{equation}
%%%%%%%%%%%%%%%%%%%%%%%%%%%%%%%%%%%%%%%%%%%%%%%%%%%%%%%%%%%%%%%%%

As is well-known or mentioned in \cite{D1}, we have
%(19)%%%%%%%%%%%%%%%%%%%%%%%%%%%%%%%%%%%%%%%%%%%%%%%%%%%%%%%%%%%%%%
\begin{equation}\label{s}
|O^-(2n,q)|=2q^{n^2-n}(q^n+1)\prod_{j=1}^{n-1}(q^{2j}-1),
\end{equation}
%%%%%%%%%%%%%%%%%%%%%%%%%%%%%%%%%%%%%%%%%%%%%%%%%%%%%%%%%%%%%%%%%
%(20)%%%%%%%%%%%%%%%%%%%%%%%%%%%%%%%%%%%%%%%%%%%%%%%%%%%%%%%%%%%%%%
\begin{equation}\label{t}
|SO^-(2n,q)|=q^{n^2-n}(q^n+1)\prod_{j=1}^{n-1}(q^{2j}-1).
\end{equation}
%%%%%%%%%%%%%%%%%%%%%%%%%%%%%%%%%%%%%%%%%%%%%%%%%%%%%%%%%%%%%%%%%

In particular, from (\ref{s}) and (\ref{t}), we have
%(21)%%%%%%%%%%%%%%%%%%%%%%%%%%%%%%%%%%%%%%%%%%%%%%%%%%%%%%%%%%%%%%
\begin{equation}\label{u}
[O^-(2n,q):SO^-(2n,q)]=2,~|SO^-(2,q)|=q+1.
\end{equation}
%%%%%%%%%%%%%%%%%%%%%%%%%%%%%%%%%%%%%%%%%%%%%%%%%%%%%%%%%%%%%%%%%

%%%%%%%%%%%%%%%%%%%%%%%%%%%%%%%%%%%%%%%%%%%%%%%%%%%%%%%%%%%%%%%%%%%%%%%%
%%%%%%%%%%%%%%%%%%%%%%%%%%%%%%%%%%%%%%%%%%%%%%%%%%%%%%%%%%%%%%%%%%%%%%%%
\section{Gauss sums for $O^{-}(2n,q)$}
%%%%%%%%%%%%%%%%%%%%%%%%%%%%%%%%%%%%%%%%%%%%%%%%%%%%%%%%%%%%%%%%%%%%%%%%
%%%%%%%%%%%%%%%%%%%%%%%%%%%%%%%%%%%%%%%%%%%%%%%%%%%%%%%%%%%%%%%%%%%%%%%%
The following notations will be employed throughout this paper.
\begin{align*}
\begin{split}
tr(x)&=x+x^3+\cdots+x^{3^{r-1}}~the~trace~function~\mathbb{F}_q\rightarrow\mathbb{F}_3,\\
 \lambda_0(x)&=e^{2\pi
ix/3}~the~canonical~additive~character~of~\mathbb{F}_3,\\
\lambda(x)&=e^{2\pi i
tr(x)/3}~the~canonical~additive~character~of~\mathbb{F}_q.
\end{split}
\end{align*}

Then any nontrivial additive character $\psi$ of $\mathbb{F}_q$ is
given by $\psi(x)=\lambda(ax)$, for a unique
$a\in\mathbb{F}_{q}^{*}.$

For any nontrivial additive character $\psi$ of $\mathbb{F}_q$ and
$a\in\mathbb{F}_{q}^{*}$, the Kloosterman sum $K_{GL(t,q)}(\psi;a)$
for $GL(t,q)$ is defined as
\begin{align*}
K_{GL(t,q)}(\psi;a)=\sum_{w\in GL(t,q)}\psi(Trw+aTrw^{-1}).
\end{align*}
Notice that, for $t=1$, $K_{GL(1,q)}(\psi;a)$ denotes the
Kloosterman sum $K(\psi;a)$.

In \cite{D3}, it is shown that $K_{GL(t,q)}(\psi;a)$ satisfies the
following recursive relation: for integers $t\geq2$,
$a\in\mathbb{F}_{q}^{*}$,
%(22)%%%%%%%%%%%%%%%%%%%%%%%%%%%%%%%%%%%%%%%%%%%%%%%%%%%%%%%%%%%%%%%
\begin{equation}\label{v}
K_{GL(t,q)}(\psi;a)=q^{t-1}K_{GL(t-1,q)}(\psi;a)K(\psi;a)+q^{2t-2}(q^{t-1}-1)K_{GL(t-2,q)}(\psi;a),
\end{equation}
%%%%%%%%%%%%%%%%%%%%%%%%%%%%%%%%%%%%%%%%%%%%%%%%%%%%%%%%%%%%%%%%%%%%
where we understand that $K_{GL(0,q)}(\psi,a)=1$.

%Proposition 2%%%%%%%%%%%%%%%%%%%%%%%%%%%%%%%%%%%%%%%%%%%%%%%%%%%%%%%
\begin{proposition}\label{B}$($\cite{D1}$)$
Let $\psi$ be a nontrivial additive character of $\mathbb{F}_q$.
For each positive integer $r$, let $\Omega_r$ be the set of all
$r\times r$ nonsingular symmetric matrices over $\mathbb{F}_q$.
Then, with $\delta_\epsilon$ as in (\ref{m}), we have
\begin{align*}
\begin{split}
b_r(\psi)&=\sum_{B\in\Omega_r}\sum_{h\in\mathbb{F}_{q}^{r \times 2}}\psi(Tr\delta_\epsilon {^t}hBh)\\
&=
\begin{cases}
q^{r(r+6)/4}\prod_{j=1}^{r/2}(q^{2j-1}-1),& \text {for r even,}\\
-q^{(r^2+4r-1)/4}\prod_{j=1}^{(r+1)/2}(q^{2j-1}-1),& \text {for r
odd.}
\end{cases}
\end{split}
\end{align*}
\end{proposition}
%%%%%%%%%%%%%%%%%%%%%%%%%%%%%%%%%%%%%%%%%%%%%%%%%%%%%%%%%%%%%%%%%%%%%

%Proposition 3%%%%%%%%%%%%%%%%%%%%%%%%%%%%%%%%%%%%%%%%%%%%%%%%%%%%%%%
\begin{proposition}\label{C}$($\cite{D2}$)$
Let $\psi$ be a nontrivial additive  character  of $\mathbb{F}_q$.
Then
\begin{enumerate}
\item \begin{align*}\label{Ca}
 \sum_{w\in SO^-(2,q)}\psi(Trw)=-K(\psi;1),
\end{align*}

\item \label{Cb}
\begin{align*}
 \sum_{w\in SO^-(2,q)}\psi(Tr\delta_1w)=q+1,
\end{align*}

\item\label{Cc}
\begin{align*}
 \sum_{w\in O^-(2,q)}\psi(Trw)=-K(\lambda;1)+q+1~(cf.(\ref{n})),
\end{align*}
where
\begin{align*}
\delta_1=\left[
           \begin{array}{cc}
             1 & 0 \\
             0 & -1 \\
           \end{array}
         \right].
\end{align*}
\end{enumerate}
\end{proposition}
%%%%%%%%%%%%%%%%%%%%%%%%%%%%%%%%%%%%%%%%%%%%%%%%%%%%%%%%%%%%%%%%%%%%%

In Section 5 of \cite{D1}, it is shown that the Gauss sum for
$SO^-(2n,q)$ is given by:
\begin{align*}
\sum_{w\in SO^-(2n,q)}\psi(Trw) \qquad \qquad \qquad \qquad \qquad
\qquad \qquad \qquad \qquad \qquad \qquad \qquad \quad \quad \quad
\end{align*}
\begin{align*}
=\sum_{\substack{0\leq r\leq n-1\\r~even}}|B_r\setminus
Q|\sum_{w\in Q}\psi(Trw\sigma_r)+\sum_{\substack{0\leq r\leq
n-1\\r~odd}}|B_r\setminus Q|\sum_{w\in Q}\psi(Tr\rho
w\sigma_r)~(cf. ~(\ref{r})) \quad \quad
\end{align*}
\begin{align*}
\begin{split}
=q^{(n-1)(n+2)/2}\Bigg{\{} \sum_{i\in SO^-(2,q)}\psi(Tri) \sum_{\substack{0\leq r\leq n-1\\r~even}}|B_r\setminus Q|q^{r(n-r-3)}b_r(\psi)K_{GL(n-1-r,q)}(\psi;1)   \\
 +\sum_{i\in SO^-(2,q)}\psi(Tr\delta_1i) \sum_{\substack{0\leq
r\leq n-1\\r~odd}}|B_r\setminus
Q|q^{r(n-r-3)}b_r(\psi)K_{GL(n-1-r,q)}(\psi;1) \Bigg{\}}
\end{split}
\end{align*}
%(23)%%%%%%%%%%%%%%%%%%%%%%%%%%%%%%%%%%%%%%%%%%%%%%%%%%%%%%%%%%%%%%%%
\begin{equation}\label{w}
\begin{split}
=-&q^{(n-1)(n+2)/2}\Bigg{\{}K(\psi;1)\sum_{\substack{0\leq r\leq
n-1\\r~even}} {[{{}_{~~~ r}^{n - 1}}]}_q
q^{rn-{\frac{1}{4}}r^2}\prod_{j=1}^{r/2}(q^{2j-1}-1)
K_{GL(n-1-r,q)}(\psi;1) \\
&+(q+1)\sum_{\substack{0\leq r\leq n-1\\r~odd}}{[{{}_{~~~ r}^{n -
1}}]}_q
q^{rn-{\frac{1}{4}}(r+1)^2}\prod_{j=1}^{(r+1)/2}(q^{2j-1}-1)
K_{GL(n-1-r,q)}(\psi;1)  \Bigg{\}}.
\end{split}
\end{equation}
%%%%%%%%%%%%%%%%%%%%%%%%%%%%%%%%%%%%%%%%%%%%%%%%%%%%%%%%%%%%%%%%%%%%%%%

Also, from Section 6 of \cite{D1} the Gauss sum for $O^-(2n,q)$ is
given by:
\begin{align*}
\sum_{w\in O^-(2n,q)}\psi(Trw) \qquad \qquad \qquad \qquad \qquad
\qquad \qquad \qquad \qquad \qquad \qquad \qquad \quad \quad \quad
\end{align*}
\begin{align*}
=\sum_{0\leq r\leq n-1}|B_r\setminus Q|\sum_{w\in
Q}\psi(Trw\sigma_r)+\sum_{0\leq r\leq n-1}|B_r\setminus
Q|\sum_{w\in Q}\psi(Tr \rho w\sigma_r)~(cf. (\ref{q})) \quad \quad
\end{align*}
\begin{align*}
=q^{(n-1)(n+2)/2}(-K(\psi;1)+q+1)\sum_{0\leq r\leq n-1}|B_r\setminus
Q|q^{r(n-r-3)}b_r(\psi)K_{GL(n-1-r,q)}(\psi;1)
\end{align*}
%(24)%%%%%%%%%%%%%%%%%%%%%%%%%%%%%%%%%%%%%%%%%%%%%%%%%%%%%%%%%%%%%%%%
\begin{equation}\label{x}
\begin{split}
=q&^{(n-1)(n+2)/2}(-K(\psi;1)+q+1)\\
&\qquad\times\Bigg{\{}\sum_{\substack{0\leq r\leq n-1\\r~even}}
{[{{}_{~~~ r}^{n - 1}}]}_q
q^{rn-{\frac{1}{4}}r^2}\prod_{j=1}^{r/2}(q^{2j-1}-1)
K_{GL(n-1-r,q)}(\psi;1) \\
&\qquad -\sum_{\substack{0\leq r\leq n-1\\r~odd}}{[{{}_{~~~ r}^{n
- 1}}]}_q
q^{rn-{\frac{1}{4}}(r+1)^2}\prod_{j=1}^{(r+1)/2}(q^{2j-1}-1)
K_{GL(n-1-r,q)}(\psi;1) \Bigg{\}}.\qquad
\end{split}
\end{equation}
%%%%%%%%%%%%%%%%%%%%%%%%%%%%%%%%%%%%%%%%%%%%%%%%%%%%%%%%%%%%%%%%%%%%%%%
Note here that the results in Proposition \ref{C} are incorporated
in (\ref{w}) and (\ref{x}).

For our purposes, we only need the following three expressions of
the Gauss sums for $SO^-(2,q)$, $O^-(2,q)$, and $SO^-(4,q)$. So we
state them separately as a theorem which follows from
(\ref{v})-(\ref{x}) by a simple change of variables. Also, for the
ease of notations, we introduce
\begin{align*}
G_1(q)=SO^-(2,q),~G_2(q)=O^-(2,q),~G_3(q)=SO^-(4,q).
\end{align*}

%Theorem (4)%%%%%%%%%%%%%%%%%%%%%%%%%%%%%%%%%%%%%%%%%%%%%%%%%%%%%%%%%%%
\begin{theorem}\label{D}
Let $\lambda$ be the canonical additive character of
$\mathbb{F}_q$, and let $a\in\mathbb{F}_{q}^{*}$. Then we have
%(25)%%%%%%%%%%%%%%%%%%%%%%%%%%%%%%%%%%%%%%%%%%%%%%%%%%%%%%%%%%%%%%%%%%
\begin{equation}\label{y}
\sum_{w\in G_1(q)}\lambda(aTrw)=-K(\lambda;a^2),
\end{equation}
%%%%%%%%%%%%%%%%%%%%%%%%%%%%%%%%%%%%%%%%%%%%%%%%%%%%%%%%%%%%%%%%%%%%%%%
%(26)%%%%%%%%%%%%%%%%%%%%%%%%%%%%%%%%%%%%%%%%%%%%%%%%%%%%%%%%%%%%%%%%%%
\begin{equation}\label{z}
\sum_{w\in G_2(q)}\lambda(aTrw)=-K(\lambda;a^2)+q+1,
\end{equation}
%%%%%%%%%%%%%%%%%%%%%%%%%%%%%%%%%%%%%%%%%%%%%%%%%%%%%%%%%%%%%%%%%%%%%%%
%(27)%%%%%%%%%%%%%%%%%%%%%%%%%%%%%%%%%%%%%%%%%%%%%%%%%%%%%%%%%%%%%%%%%%
\begin{equation}\label{a1}
\sum_{w\in G_3(q)}\lambda(aTrw)=-q^2(K(\lambda;a^2)^2+q^3-q).
\end{equation}
%%%%%%%%%%%%%%%%%%%%%%%%%%%%%%%%%%%%%%%%%%%%%%%%%%%%%%%%%%%%%%%%%%%%%%%
\end{theorem}
%%%%%%%%%%%%%%%%%%%%%%%%%%%%%%%%%%%%%%%%%%%%%%%%%%%%%%%%%%%%%%%%%%%%%%%

%Proposition (5)%%%%%%%%%%%%%%%%%%%%%%%%%%%%%%%%%%%%%%%%%%%%%%%%%%%%%%%%%%%
\begin{proposition}\label{E}$[$$7$, $( 5.3-5 ) ]$
Let $\lambda$ be the canonical additive character of
$\mathbb{F}_q$, $m\in \mathbb{Z}_{\geq0}$, $\beta\in\mathbb{F}_q$.
Then
%(28)%%%%%%%%%%%%%%%%%%%%%%%%%%%%%%%%%%%%%%%%%%%%%%%%%%%%%%%%%%%%%%%%%%
\begin{equation}\label{b1}
\sum_{a\in\mathbb{F}_{q}^{*}}\lambda(-a\beta)K(\lambda;a^2)^m=q\delta(m,q;\beta)-(q-1)^m,
\end{equation}
%%%%%%%%%%%%%%%%%%%%%%%%%%%%%%%%%%%%%%%%%%%%%%%%%%%%%%%%%%%%%%%%%%%%%%%
where, for $m\geq1$,
%(29)%%%%%%%%%%%%%%%%%%%%%%%%%%%%%%%%%%%%%%%%%%%%%%%%%%%%%%%%%%%%%%%%%%
\begin{equation}\label{c1}
\delta(m,q;\beta)=|\{(\alpha_{1}^{},\cdots,\alpha_{m}^{})\in(\mathbb{F}_{q}^{*})^m|\alpha_{1}^{}+\alpha_{1}^{-1}+
\cdots,\alpha_{m}^{}+\alpha_{m}^{-1}=\beta\}|,
\end{equation}
%%%%%%%%%%%%%%%%%%%%%%%%%%%%%%%%%%%%%%%%%%%%%%%%%%%%%%%%%%%%%%%%%%%%%%%
and
%(30)%%%%%%%%%%%%%%%%%%%%%%%%%%%%%%%%%%%%%%%%%%%%%%%%%%%%%%%%%%%%%%%%%%
\begin{equation}\label{d1}
\delta(0,q;\beta)=
\begin{cases}
1,& \text {$\beta$=0,}\\
0,& \text {otherwise.}
\end{cases}
\end{equation}
%%%%%%%%%%%%%%%%%%%%%%%%%%%%%%%%%%%%%%%%%%%%%%%%%%%%%%%%%%%%%%%%%%%%%%%
\end{proposition}
%%%%%%%%%%%%%%%%%%%%%%%%%%%%%%%%%%%%%%%%%%%%%%%%%%%%%%%%%%%%%%%%%%%%%%%

%Remark 6%%%%%%%%%%%%%%%%%%%%%%%%%%%%%%%%%%%%%%%%%%%%%%%%%%%%%%%%%%%%%%%
\begin{remark}\label{F}
Here one notes that
%(31)%%%%%%%%%%%%%%%%%%%%%%%%%%%%%%%%%%%%%%%%%%%%%%%%%%%%%%%%%%%%%%%%%%
\begin{equation}\label{e1}
\begin{split}
\delta(1,q;\beta)&=|\{x\in\mathbb{F}_{q}^{}|x^2-\beta x+1=0\}|\\
&=
\begin{cases}
2,& \text {if $\beta^2-1\neq$$0$ is a square,}\\
1,& \text {if $\beta^2-1$= $0$,}\\
0,& \text {if $\beta^2-1$ is a nonsquare.}
\end{cases}
\end{split}
\end{equation}
%%%%%%%%%%%%%%%%%%%%%%%%%%%%%%%%%%%%%%%%%%%%%%%%%%%%%%%%%%%%%%%%%%%%%%%
\end{remark}
%%%%%%%%%%%%%%%%%%%%%%%%%%%%%%%%%%%%%%%%%%%%%%%%%%%%%%%%%%%%%%%%%%%%%%%%

Let $G(q)$ be one of finite classical groups over $\mathbb{F}_q$.
Then we put, for each $\beta\in\mathbb{F}_q$,
\begin{align*}
N_{G(q)}(\beta)=|\{w\in G(q)|~Tr(w)=\beta\}|.
\end{align*}
Then it is easy to see that
%(32)%%%%%%%%%%%%%%%%%%%%%%%%%%%%%%%%%%%%%%%%%%%%%%%%%%%%%%%%%%%%%%%%%%%
\begin{equation}\label{f1}
qN_{G(q)}(\beta)=|G(q)|+\sum_{a\in\mathbb{F}_{q}^{*}}\lambda(-a\beta)
\sum_{w\in G(q)}\lambda(aTrw).
\end{equation}
%%%%%%%%%%%%%%%%%%%%%%%%%%%%%%%%%%%%%%%%%%%%%%%%%%%%%%%%%%%%%%%%%%%%%%%%
For brevity, we write
%(33)%%%%%%%%%%%%%%%%%%%%%%%%%%%%%%%%%%%%%%%%%%%%%%%%%%%%%%%%%%%%%%%%%%%
\begin{equation}\label{g1}
n_1(\beta)=N_{G_1(q)}(\beta),~n_2(\beta)=N_{G_2(q)}(\beta),~n_3(\beta)=N_{G_3(q)}(\beta).
\end{equation}
%%%%%%%%%%%%%%%%%%%%%%%%%%%%%%%%%%%%%%%%%%%%%%%%%%%%%%%%%%%%%%%%%%%%%%%%
Using (\ref{f1}), (\ref{m1}), (\ref{y})-(\ref{a1}), one derives the
following.

%Proposition 7%%%%%%%%%%%%%%%%%%%%%%%%%%%%%%%%%%%%%%%%%%%%%%%%%%%%%%%%%%
\begin{proposition}\label{G}
Let $n_1(\beta)$, $n_2(\beta)$, $n_3(\beta)$ be as in (\ref{g1}).
Then we have the following.
\begin{enumerate}
\item\label{Ga}
%(34)%%%%%%%%%%%%%%%%%%%%%%%%%%%%%%%%%%%%%%%%%%%%%%%%%%%%%%%%%%%%%%%%%%%
\begin{equation}\label{h1}
n_1(\beta)=2-\delta(1,q;\beta)=
\begin{cases}
0,& \text {if $\beta^2-1\neq$$0$ is a square,}\\
1,& \text {if $\beta^2-1$=$0$,}\\
2,& \text {if $\beta^2-1$ is a nonsquare.}
\end{cases}
\end{equation}
%%%%%%%%%%%%%%%%%%%%%%%%%%%%%%%%%%%%%%%%%%%%%%%%%%%%%%%%%%%%%%%%%%%%%%%%

\item\label{Gb}
\begin{enumerate}
\item\label{Gba}
For $r$ is even,
%(35)%%%%%%%%%%%%%%%%%%%%%%%%%%%%%%%%%%%%%%%%%%%%%%%%%%%%%%%%%%%%%%%%%%%
\begin{equation}\label{i1}
n_2(\beta)=
\begin{cases}
0,& \text {if $\beta^2-1\neq-1,0$ is a square,}\\
q+1,& \text {if $\beta^2-1=-1$,}\\
1,& \text {if $\beta^2-1=0$,}\\
2,& \text {if $\beta^2-1$ is a nonsquare.}
\end{cases}
\end{equation}
%%%%%%%%%%%%%%%%%%%%%%%%%%%%%%%%%%%%%%%%%%%%%%%%%%%%%%%%%%%%%%%%%%%%%%%%

\item\label{Gbb}
For $r$ is odd,
%(36)%%%%%%%%%%%%%%%%%%%%%%%%%%%%%%%%%%%%%%%%%%%%%%%%%%%%%%%%%%%%%%%%%%%
\begin{equation}\label{j1}
n_2(\beta)=
\begin{cases}
0,& \text {if $\beta^2-1\neq0$ is a square,}\\
1,& \text {if $\beta^2-1=0$,}\\
2,& \text {if $\beta^2-1\neq-1$ is a nonsquare,}\\
q+3,& \text {if $\beta^2-1=-1$.}
\end{cases}
\end{equation}
%%%%%%%%%%%%%%%%%%%%%%%%%%%%%%%%%%%%%%%%%%%%%%%%%%%%%%%%%%%%%%%%%%%%%%%%
\end{enumerate}

\item\label{Gc}
%(37)%%%%%%%%%%%%%%%%%%%%%%%%%%%%%%%%%%%%%%%%%%%%%%%%%%%%%%%%%%%%%%%%%%%
\begin{equation}\label{k1}
n_3(\beta)=
\begin{cases}
-q^2\delta(2,q;\beta)+q^5+q^4+q^3-3q^2,& \text {if $\beta\neq0$,}\\
-q^2\delta(2,q;0)+q^4+2q^3-3q^2,& \text {if $\beta=0$.}
\end{cases}
\end{equation}
%%%%%%%%%%%%%%%%%%%%%%%%%%%%%%%%%%%%%%%%%%%%%%%%%%%%%%%%%%%%%%%%%%%%%%%%
\end{enumerate}
Here
$\delta(2,q;\beta)=|\{(\alpha_{1}^{},\alpha_{2}^{})\in(\mathbb{F}_{q}^{*})^2|\alpha_{1}^{}+
\alpha_{1}^{-1}+\alpha_{2}^{}+\alpha_{2}^{-1}=\beta\}|$.
\end{proposition}

\begin{proof}
Here we only provide the proof for (\ref{Gb}). The others can be
shown in a similar manner. From (\ref{z}), (\ref{b1}), (\ref{f1})
and (\ref{m1}), we have
\begin{align*}
qn_2(\beta)=-q\delta(1,q;\beta)+3q+1+(q+1)\sum_{a\in\mathbb{F}_{q}^{*}}
\lambda(-a\beta),
\end{align*}
and hence
%(38)%%%%%%%%%%%%%%%%%%%%%%%%%%%%%%%%%%%%%%%%%%%%%%%%%%%%%%%%%%%%%%%%%%%
\begin{equation}\label{l1}
n_2(\beta)=
\begin{cases}
-\delta(1,q;\beta)+2,& \text {if $\beta\neq0$,}\\
-\delta(1,q;0)+q+3,& \text {if $\beta=0$.}
\end{cases}
\end{equation}
%%%%%%%%%%%%%%%%%%%%%%%%%%%%%%%%%%%%%%%%%%%%%%%%%%%%%%%%%%%%%%%%%%%%%%%%
Noting now that, for $\beta=0$, $\beta^2-1=-1$ is a square in
$\mathbb{F}_{q}^{}=\mathbb{F}_{3^r}^{}$ if and only if $r$ is even,
the results now follow from Remark \ref{F}.
\end{proof}
%%%%%%%%%%%%%%%%%%%%%%%%%%%%%%%%%%%%%%%%%%%%%%%%%%%%%%%%%%%%%%%%%%%%%%%%

%%%%%%%%%%%%%%%%%%%%%%%%%%%%%%%%%%%%%%%%%%%%%%%%%%%%%%%%%%%%%%%%%%%%%%%%
%%%%%%%%%%%%%%%%%%%%%%%%%%%%%%%%%%%%%%%%%%%%%%%%%%%%%%%%%%%%%%%%%%%%%%%%
\section{Construction of codes}
%%%%%%%%%%%%%%%%%%%%%%%%%%%%%%%%%%%%%%%%%%%%%%%%%%%%%%%%%%%%%%%%%%%%%%%%
%%%%%%%%%%%%%%%%%%%%%%%%%%%%%%%%%%%%%%%%%%%%%%%%%%%%%%%%%%%%%%%%%%%%%%%%
Let
%(39)%%%%%%%%%%%%%%%%%%%%%%%%%%%%%%%%%%%%%%%%%%%%%%%%%%%%%%%%%%%%%%%%%%%
\begin{equation}\label{m1}
N_1=|G_1(q)|=q+1,~N_2=|G_2(q)|=2(q+1),~N_3=|G_3(q)|=q^2(q^4-1).
\end{equation}
%%%%%%%%%%%%%%%%%%%%%%%%%%%%%%%%%%%%%%%%%%%%%%%%%%%%%%%%%%%%%%%%%%%%%%%%
Here we will construct three ternary linear codes $C(G_1(q))$ of
length $N_1$, $C(G_2(q))$ of length $N_2$, and $C(G_3(q))$ of length
$N_3$, respectively associated with the orthogonal groups $G_1(q)$,
$G_2(q)$, and $G_3(q)$.

By abuse of notations, for $i=1,2,3$, let $g_1,g_2,\cdots,g_{N_i}$
be a fixed ordering of the elements in the group $G_i(q)$. \\
Also, for $i=1,2,3$, we put
\begin{align*}
v_i=(Trg_1,Trg_2,\cdots,Trg_{N_i})\in\mathbb{F}_{q}^{N_i}.
\end{align*}
Then, for $i=1,2,3$, the ternary linear code $C(G_i(q))$ is defined
as
%(40)%%%%%%%%%%%%%%%%%%%%%%%%%%%%%%%%%%%%%%%%%%%%%%%%%%%%%%%%%%%%%%%%%%%
\begin{equation}\label{n1}
C(G_i(q))=\{u\in\mathbb{F}_{2}^{N_i}|u\cdot v_i=0\},
\end{equation}
%%%%%%%%%%%%%%%%%%%%%%%%%%%%%%%%%%%%%%%%%%%%%%%%%%%%%%%%%%%%%%%%%%%%%%%%
where the dot denotes the usual inner product in
$\mathbb{F}_{q}^{N_i}$.

The following Delsarte's theorem is well-known.

%Theorem 8%%%%%%%%%%%%%%%%%%%%%%%%%%%%%%%%%%%%%%%%%%%%%%%%%%%%%%%%%%%%%%%
\begin{theorem}$([13])$\label{H}
Let $B$ be a linear code over $\mathbb{F}_q$. Then
\begin{align*}
(B|_{\mathbb{F}_{3}})^\bot=tr(B^\bot).
\end{align*}
In view of this theorem, the dual $C(G_i(q))^\bot ~(i=1,2,3)$ is
given by
\begin{equation}\label{o1}
C(G_i(q))^\bot=\{c_i(a)=(tr(aTrg_1),\cdots,tr(aTrg_{N_i}))|a\in\mathbb{F}_q\}.
\end{equation}
\end{theorem}
%%%%%%%%%%%%%%%%%%%%%%%%%%%%%%%%%%%%%%%%%%%%%%%%%%%%%%%%%%%%%%%%%%%%%%%%%

%Lemma 8%%%%%%%%%%%%%%%%%%%%%%%%%%%%%%%%%%%%%%%%%%%%%%%%%%%%%%%%%%%%%%%%%
\begin{lemma}\label{I}
Let $\delta(1,q;\beta)$ be as in (\ref{c1}) and (\ref{e1}), and let
$a\in\mathbb{F}_{q}^{*}$. Then we have
%(42)%%%%%%%%%%%%%%%%%%%%%%%%%%%%%%%%%%%%%%%%%%%%%%%%%%%%%%%%%%%%%%%%%%%%%
\begin{equation}\label{p1}
\sum_{\beta\in\mathbb{F}_{q}}\delta(1,q;\beta)\lambda(a\beta)=K(\lambda;a^2).
\end{equation}
\end{lemma}
%%%%%%%%%%%%%%%%%%%%%%%%%%%%%%%%%%%%%%%%%%%%%%%%%%%%%%%%%%%%%%%%%%%%%%%%%%
\begin{proof}
The LHS of (\ref{p1}) is equal to
\begin{align*}
\begin{split}
&\sum_{\beta\in\mathbb{F}_q}(q^{-1}\sum_{x\in\mathbb{F}_q}\sum_{\alpha\in\mathbb{F}_q}\lambda(\alpha(x^2-\beta
x+1)))\lambda(a\beta)\\
&=q^{-1}\sum_{x\in\mathbb{F}_q}\sum_{\alpha\in\mathbb{F}_q}\lambda(\alpha(x^2+1))\sum_{\beta\in\mathbb{F}_q}\lambda(\beta(a-\alpha x))\\
&=\sum_{\alpha\in\mathbb{F}_{q}^{*}}\lambda(\alpha(\alpha^{-2}a^2+1))\\
&=K(\lambda;a^2).
\end{split}
\end{align*}
\end{proof}
%%%%%%%%%%%%%%%%%%%%%%%%%%%%%%%%%%%%%%%%%%%%%%%%%%%%%%%%%%%%%%%%%%%%%%%%%

%Proposition (%%%%%%%%%%%%%%%%%%%%%%%%%%%%%%%%%%%%%%%%%%%%%%%%%%%%%%%%%%%%
\begin{proposition}\label{J}
The map $\mathbb{F}_q \rightarrow C(G_i(q))^\bot~(a\mapsto c_i(a))$
is an $\mathbb{F}_3$-linear isomorphism, for every $q=3^r$, and
$i=1,2,3$.
\end{proposition}
\begin{proof}
As the proofs for $i=1$ and $i=2$ are similar, we will treat only
the cases $i=2$ and $i=3$. Let $i=2$. The map is clearly
$\mathbb{F}_3$-linear and surjective. Let $a$ be in the kernel of
the map. Then, $tr(aTrg)=0$, for all $g\in G_2(q)$. Assume that
$a\neq0$. Then
\begin{align*}
\begin{split}
2(q+1)=|G_2(q)|&=\sum_{g\in G_2(q)}e^{2\pi itr(aTrg)/3}\\
&=\sum_{\beta\in\mathbb{F}_q}n_2(\beta)\lambda(a\beta)\\
&=q+1+\sum_{\beta\in\mathbb{F}_q}(-\delta(1,q;\beta)+2)\lambda(a\beta)~(from~(\ref{l1}))\\
&=q+1-\sum_{\beta\in\mathbb{F}_q}\delta(1,q;\beta)\lambda(a\beta)\\
&=q+1-K(\lambda;a^2).~(by~(\ref{p1}))
\end{split}
\end{align*}
So $q+1=-K(\lambda;a^2)$. This yields from Weil bound (\ref{a})
that $q+1\leq 2\sqrt{q}$, equivalently $q=1$. As $q=3^r\geq3$, we
must have $a=0$. Now, let $i=3$. Then again the map is
$\mathbb{F}_3$-linear and surjective. Let $a$ be in the kernal of
the map. Then $tr(aTrg)=0$, for all $g\in G_3(q)$. Since
$n_3(\beta)>0$, for all $\beta\in\mathbb{F}_{q}^{}$(cf.
(\ref{k1})), this in turn yields that $tr(a\beta)=0$. As the trace
function $\mathbb{F}_q\rightarrow\mathbb{F}_3$ is surjective, we
must have $a=0$.
\end{proof}
%%%%%%%%%%%%%%%%%%%%%%%%%%%%%%%%%%%%%%%%%%%%%%%%%%%%%%%%%%%%%%%%%%%%%%%%%%

%%%%%%%%%%%%%%%%%%%%%%%%%%%%%%%%%%%%%%%%%%%%%%%%%%%%%%%%%%%%%%%%%%%%%%%%
%%%%%%%%%%%%%%%%%%%%%%%%%%%%%%%%%%%%%%%%%%%%%%%%%%%%%%%%%%%%%%%%%%%%%%%%
\section{Recusive formulas for power moments of Kloosterman sums with square arguments}
%%%%%%%%%%%%%%%%%%%%%%%%%%%%%%%%%%%%%%%%%%%%%%%%%%%%%%%%%%%%%%%%%%%%%%%%
%%%%%%%%%%%%%%%%%%%%%%%%%%%%%%%%%%%%%%%%%%%%%%%%%%%%%%%%%%%%%%%%%%%%%%%%
 In this section, we will be able to find, via Pless power moment identity,
 recursive formulas for the power moments of Kloosterman sums with
square arguments and even power moments of those with square
arguments  in terms of the frequencies of weights in $C(G_i(q))$,
for each $i=1,2,3$.

%Theorem 11%%%%%%%%%%%%%%%%%%%%%%%%%%%%%%%%%%%%%%%%%%%%%%%%%%%%%%%%%%%%%%%%%%
\begin{theorem}\label{K}
$($Pless power moment identity, \cite{MS1}$)$
 Let $B$ be an $q$-ary $[n,k]$ code, and let $B_i$(resp. $B_i^\bot$)
 denote the number of codewords of weight $i$ in $B$(resp. in
 $B^\bot$). Then, for $h=0,1,2,\cdots,$
%(43)%%%%%%%%%%%%%%%%%%%%%%%%%%%%%%%%%%%%%%%%%%%%%%%%%%%%%%%%%%%%%%%%%%%%%%%%
\begin{equation}\label{q1}
\sum_{j=0}^{n}j^hB_j=\sum_{j=0}^{min\{n,h\}}(-1)^jB_j^\bot
\sum_{t=j}^{h}t!S(h,t)q^{k-t}(q-1)^{t-j}{\binom{n-j}{n-t}},
\end{equation}
%%%%%%%%%%%%%%%%%%%%%%%%%%%%%%%%%%%%%%%%%%%%%%%%%%%%%%%%%%%%%%%%%%%%%%%%%%%%%
where $S(h,t)$ is the Stirling number of the second kind defined in
(\ref{g}).
\end{theorem}
%%%%%%%%%%%%%%%%%%%%%%%%%%%%%%%%%%%%%%%%%%%%%%%%%%%%%%%%%%%%%%%%%%%%%%%%%%%%%%

%Lemma 12%%%%%%%%%%%%%%%%%%%%%%%%%%%%%%%%%%%%%%%%%%%%%%%%%%%%%%%%%%%%%%%%%%%%%
\begin{lemma}\label{L}
Let $c_i(a)=(tr(aTrg_1),\cdots,tr(aTrg_{N_i}))\in C(G_i(q))^\bot$,
for $a\in\mathbb{F}_{q}^{*}$, and $i=1,2,3$. Then the Hamming
weight $w(c_i(a))$ can be expressed as follows:
\begin{enumerate}
\item\label{La}
%(44)%%%%%%%%%%%%%%%%%%%%%%%%%%%%%%%%%%%%%%%%%%%%%%%%%%%%%%%%%%%%%%%%%%%%%%%%
\begin{equation}\label{r1}
w(c_i(a))={\frac{2}{3}}(q+1+K(\lambda;a^2)),~for~i=1,2,
\end{equation}
%%%%%%%%%%%%%%%%%%%%%%%%%%%%%%%%%%%%%%%%%%%%%%%%%%%%%%%%%%%%%%%%%%%%%%%%%%%%%
\item\label{Lb}
%(45)%%%%%%%%%%%%%%%%%%%%%%%%%%%%%%%%%%%%%%%%%%%%%%%%%%%%%%%%%%%%%%%%%%%%%%%%
\begin{equation}\label{s1}
w(c_3(a))={\frac{2}{3}}q^2(K(\lambda;a^2)^2+q^4+q^3-q-1).
\end{equation}
%%%%%%%%%%%%%%%%%%%%%%%%%%%%%%%%%%%%%%%%%%%%%%%%%%%%%%%%%%%%%%%%%%%%%%%%%%%%%
\end{enumerate}
\begin{proof} For $i=1,2,3,$
\begin{align*}
\begin{split}
w(c_i(a))&=\sum_{j=1}^{N_i}(1-{\frac{1}{3}}\sum_{\alpha\in\mathbb{F}_3}\lambda_0(\alpha tr(aTrg_j)))\\
&=N_i-{\frac{1}{3}}\sum_{\alpha\in\mathbb{F}_3}\sum_{w\in
G_i(q)}\lambda(\alpha aTrw)\\
&={\frac{2}{3}}N_i-{\frac{1}{3}}\sum_{\alpha\in\mathbb{F}_3^*}\sum_{w\in
G_i(q)}\lambda(\alpha aTrw).
\end{split}
\end{align*}
Our results now follow from (\ref{m1}) and (\ref{y})-(\ref{a1}).
\end{proof}
\end{lemma}
%%%%%%%%%%%%%%%%%%%%%%%%%%%%%%%%%%%%%%%%%%%%%%%%%%%%%%%%%%%%%%%%%%%%%%%%%%%%%%%%

Fix $i(i=1,2,3)$, and let
$u=(u_1,\cdots,u_{N_i})\in\mathbb{F}_3^{N_i}$, with $\nu_\beta$ 1's
and $\mu_\beta$ 2's in the coordinate places where $Tr(g_j)=\beta$,
for each $\beta\in\mathbb{F}_{q}^{}$. Then we see from the
definition of the code $C(G_i(q))$(cf. (\ref{n1})) that $u$ is a
codeword with weight $j$ if and only if
$\sum_{\beta\in\mathbb{F}_q}\nu_\beta+\sum_{\beta\in\mathbb{F}_q}\mu_\beta=j$
and $\sum_{\beta\in\mathbb{F}_q}\nu_\beta
\beta=\sum_{\beta\in\mathbb{F}_q}\mu_\beta \beta$(an identity in
$\mathbb{F}_q$). Note that there are
$\prod_{\beta\in\mathbb{F}_q}{\binom{n_i(\beta)}{\nu_\beta,\mu_\beta}}$(cf.
(\ref{c}), (\ref{d})) many such codewords with weight $j$. Now, we
get the following formulas in (\ref{t1})-(\ref{w1}), by using the
explicit values of $n_i(\beta)$ in (\ref{h1})-(\ref{k1}).

%Theorem 13%%%%%%%%%%%%%%%%%%%%%%%%%%%%%%%%%%%%%%%%%%%%%%%%%%%%%%%%%%%%%%%%%%%%%%
\begin{theorem}\label{M}
Let $\{C_{i,j}\}_{j=0}^{N_i}$ be the weight distribution of
$C(G_i(q))$, for $i=1,2,3$. Then we have the following.

\item\label{Ma}
%(46)%%%%%%%%%%%%%%%%%%%%%%%%%%%%%%%%%%%%%%%%%%%%%%%%%%%%%%%%%%%%%%%%%%%%%%%%%%%%%
\begin{equation}\label{t1}
C_{1,j}=\sum{\binom{1}{\nu_1,\mu_1}}{\binom{1}{\nu_{-1},\mu_{-1}}}\prod_{\beta^2-1~nonsquare}
{\binom{2}{\nu_\beta,\mu_\beta}}~(j=0,\cdots,N_1),
\end{equation}
%%%%%%%%%%%%%%%%%%%%%%%%%%%%%%%%%%%%%%%%%%%%%%%%%%%%%%%%%%%%%%%%%%%%%%%%%%%%%%%%%%
where the sum is over all the sets of nonnegative integers
$\{\nu_\beta\}_{\beta^2-1~nonsquare}\cup\{\nu_{\pm1}\}$ and
$\{\mu_\beta\}_{\beta^2-1~nonsquare}\cup\{\mu_{\pm1}\}$ satisfying
\begin{align*}
\nu_1+\nu_{-1}+\sum_{\beta^2-1~nonsquare}\nu_\beta+\mu_1+\mu_{-1}+\sum_{\beta^2-1~nonsquare}\mu_\beta=j
\end{align*}
and
\begin{align*}
\nu_1 1+\nu_{-1}(-1)+\sum_{\beta^2-1~nonsquare}\nu_\beta
\beta=\mu_1 1+\mu_{-1}(-1)+\sum_{\beta^2-1~nonsquare}\mu_\beta
\beta.
\end{align*}

\item\label{Mb}

\item\label{Mba}
For $r$ even,
%(47)%%%%%%%%%%%%%%%%%%%%%%%%%%%%%%%%%%%%%%%%%%%%%%%%%%%%%%%%%%%%%%%%%%%%%%%%%%%%%
\begin{equation}\label{u1}
C_{2,j}=\sum{\binom{q+1}{\nu_0,\mu_0}}{\binom{1}{\nu_1,\mu_1}}{\binom{1}{\nu_{-1},\mu_{-1}}}\prod_{\beta^2-1~nonsquare}
{\binom{2}{\nu_\beta,\mu_\beta}}~(j=0,\cdots,N_1),
\end{equation}
%%%%%%%%%%%%%%%%%%%%%%%%%%%%%%%%%%%%%%%%%%%%%%%%%%%%%%%%%%%%%%%%%%%%%%%%%%%%%%%%%%
 where the sum is over all the sets of nonnegative integers
$\{\nu_\beta\}_{\beta^2-1~nonsquare}$\\$\cup\{\nu_{\pm1}\}\cup\{\nu_0\}$
and
$\{\mu_\beta\}_{\beta^2-1~nonsquare}\cup\{\mu_{\pm1}\}\cup\{\mu_0\}$
satisfying
\begin{align*}
\begin{split}
\nu_0+\nu_1+\nu_{-1}+&\sum_{\beta^2-1~nonsquare}\nu_\beta\\
&+\mu_0+\mu_1+\mu_{-1}+\sum_{\beta^2-1~nonsquare}\mu_\beta=j
\end{split}
\end{align*}
and
\begin{align*}
\begin{split}
\nu_1 1+\nu_{-1}(-1)+&\sum_{\beta^2-1~nonsquare}\nu_\beta \beta\\
&=\mu_1 1+\mu_{-1}(-1)+\sum_{\beta^2-1~nonsquare}\mu_\beta \beta.
\end{split}
\end{align*}

\item\label{Mbb}
For $r$ odd,
%(48)%%%%%%%%%%%%%%%%%%%%%%%%%%%%%%%%%%%%%%%%%%%%%%%%%%%%%%%%%%%%%%%%%%%%%%%%%%%%%
\begin{equation}\label{v1}
C_{2,j}=\sum{\binom{1}{\nu_1,\mu_1}}{\binom{1}{\nu_{-1},\mu_{-1}}}\prod_{\beta^2-1\neq-1~nonsquare}
{\binom{2}{\nu_\beta,\mu_\beta}}{\binom{q+3}{\nu_0,\mu_0}}~(j=0,\cdots,N_2),
\end{equation}
%%%%%%%%%%%%%%%%%%%%%%%%%%%%%%%%%%%%%%%%%%%%%%%%%%%%%%%%%%%%%%%%%%%%%%%%%%%%%%%%%%
where the sum is over all the sets of nonnegative integers
$\{\nu_\beta\}_{\beta^2-1\neq-1~nonsquare}$\\$\cup\{\nu_{\pm1}\}\cup\{\nu_0\}$
and
$\{\mu_\beta\}_{\beta^2-1\neq-1~nonsquare}\cup\{\mu_{\pm1}\}\cup\{\mu_0\}$
satisfying
\begin{align*}
\begin{split}
\nu_0+\nu_1+\nu_{-1}+&\sum_{\beta^2-1\neq-1~nonsquare}\nu_\beta\\
&+\mu_0+\mu_1+\mu_{-1}+\sum_{\beta^2-1\neq-1~nonsquare}\mu_\beta=j
\end{split}
\end{align*}
and
\begin{align*}
\begin{split}
\nu_1 1+\nu_{-1}(-1)+&\sum_{\beta^2-1\neq-1~nonsquare}\nu_\beta \beta\\
&=\mu_1 1+\mu_{-1}(-1)+\sum_{\beta^2-1\neq-1~nonsquare}\mu_\beta
\beta.
\end{split}
\end{align*}

\item\label{Mc}
%(49)%%%%%%%%%%%%%%%%%%%%%%%%%%%%%%%%%%%%%%%%%%%%%%%%%%%%%%%%%%%%%%%%%%%%%%%%%%%%%
\begin{equation}\label{w1}
\begin{split}
C_{3,j}=\sum&{\binom{-q^2\delta(2,q;0)+q^4+2q^3-3q^2}{\nu_0,\mu_0}}\\
&\times\prod_{\beta\in\mathbb{F}_{q}^{*}}
{\binom{-q^2\delta(2,q;\beta)+q^5+q^4+q^3-3q^2}{\nu_\beta,\mu_\beta}}~(j=0,\cdots,N_3),
\end{split}
\end{equation}
%%%%%%%%%%%%%%%%%%%%%%%%%%%%%%%%%%%%%%%%%%%%%%%%%%%%%%%%%%%%%%%%%%%%%%%%%%%%%%%%%%
where the sum runs over all the sets of nonnegative integers
$\{\nu_\beta\}_{\beta\in\mathbb{F}_q}$ and
$\{\mu_\beta\}_{\beta\in\mathbb{F}_q}$ satisfying
\begin{align*}
\sum_{\beta\in\mathbb{F}_q} \nu_\beta
+\sum_{\beta\in\mathbb{F}_q}\mu_\beta=j~and~
\sum_{\beta\in\mathbb{F}_q} \nu_\beta
\beta=\sum_{\beta\in\mathbb{F}_q}\mu_\beta \beta,
\end{align*}
and, for every $\beta\in\mathbb{F}_q$,
\begin{align*}
\delta(2,q;\beta)=|\{(\alpha_{1}^{},\alpha_{2}^{})\in(\mathbb{F}_{q}^{*})^2|
\alpha_{1}^{}+\alpha_{1}^{-1}+\alpha_{2}^{}+\alpha_{2}^{-1}=\beta\}|.
\end{align*}

\end{theorem}
%%%%%%%%%%%%%%%%%%%%%%%%%%%%%%%%%%%%%%%%%%%%%%%%%%%%%%%%%%%%%%%%%%%%%%%%%%%%%%%%%

We now apply the Pless power moment identity in (\ref{q1}) to each
$C(G_i(q))^\bot$, for $i=1,2,3$, in order to obtain the results in
Theorem \ref{Aa} about recursive formulas. Then the left hand side
of the identity in (\ref{q1}) is equal to
%(50)%%%%%%%%%%%%%%%%%%%%%%%%%%%%%%%%%%%%%%%%%%%%%%%%%%%%%%%%%%%%%%%%%%%%%%%%%%%%
\begin{equation}\label{x1}
\sum_{a\in\mathbb{F_{q}^{*}}}w(c_i(a))^h,
\end{equation}
%%%%%%%%%%%%%%%%%%%%%%%%%%%%%%%%%%%%%%%%%%%%%%%%%%%%%%%%%%%%%%%%%%%%%%%%%%%%%%%%%
with the $w(c_i(a))$ in each case given by (\ref{r1}) and
(\ref{s1}).

For $i=1,2$, (\ref{x1}) is
\begin{align*}
\begin{split}
&(\frac{2}{3})^h\sum_{a\in\mathbb{F}_{q}^{*}}(q+1+K(\lambda;a^2))^h\\
&=(\frac{2}{3})^h\sum_{a\in\mathbb{F}_{q}^{*}}\sum_{j=0}^{h}{\binom{h}{j}}(q+1)^{h-j}K(\lambda;a^2)^j
\end{split}
\end{align*}
%(51)%%%%%%%%%%%%%%%%%%%%%%%%%%%%%%%%%%%%%%%%%%%%%%%%%%%%%%%%%%%%%%%%%%%%%%%%%%%%
\begin{equation}\label{y1}
=2(\frac{2}{3})^h\sum_{j=0}^{h}{\binom{h}{j}}(q+1)^{h-j}SK^j.\qquad
\quad
\end{equation}
%%%%%%%%%%%%%%%%%%%%%%%%%%%%%%%%%%%%%%%%%%%%%%%%%%%%%%%%%%%%%%%%%%%%%%%%%%%%%%%%%
Similarly, for $i=3$, (\ref{x1}) equals
%(52)%%%%%%%%%%%%%%%%%%%%%%%%%%%%%%%%%%%%%%%%%%%%%%%%%%%%%%%%%%%%%%%%%%%%%%%%%%%%
\begin{equation}\label{z1}
=2(\frac{2q^2}{3})^h\sum_{j=0}^{h}{\binom{h}{j}}(q^4+q^3-q-1)^{h-j}SK^{2j}.
\end{equation}
%%%%%%%%%%%%%%%%%%%%%%%%%%%%%%%%%%%%%%%%%%%%%%%%%%%%%%%%%%%%%%%%%%%%%%%%%%%%%%%%%
Here one has to separate the term corresponding to $j=h$ in
(\ref{y1}) and (\ref{z1}), and note
$dim_{\mathbb{F}_3}C(G_i(q))^\bot=r$.

\bibliographystyle{amsplain}

\end{document}